\documentclass{article}
\usepackage{arxiv}

\usepackage{amsmath}
    \numberwithin{equation}{section}
\usepackage{hyperref}
\usepackage[utf8]{inputenc} %
\usepackage[T1]{fontenc}    %
\usepackage{booktabs}       %
\usepackage{amsfonts}       %
\usepackage{nicefrac}       %
\usepackage{microtype}      %
\usepackage{mathabx}
\usepackage{bm}
\usepackage{stmaryrd}
\usepackage{xspace}
\usepackage{arydshln}  %
\usepackage{graphicx}
\usepackage{adjustbox}
\usepackage{algorithm}     %
\usepackage{algpseudocode}
\usepackage{cleveref}      %
\usepackage[title]{appendix}
\usepackage{graphicx,wrapfig,lipsum}
\usepackage{subcaption}
\usepackage{cite}

\usepackage{printlen} %

\usepackage{xcolor}

\graphicspath{ {.} }

\newcommand{\norm}[1]{\lVert#1\rVert}

\DeclareMathOperator{\ujump}{\llbracket \disp \rrbracket}
\DeclareMathOperator{\ujumpn}{\llbracket {u}_n \rrbracket}
\newcommand{\ujumpt}[0]{\llbracket \dot{\disp}_\tau \rrbracket}
\DeclareMathOperator{\Lame}{\lambda_\text{\textnormal{\LameName}}}

\newcommand{\perm}[0]{\mathbf{K}}
\newcommand{\visc}[0]{\mu^{f}}
\DeclareMathOperator{\lambdat}{\bm{\lambda}_\tau}
\DeclareMathOperator{\lambdan}{{\lambda}_n}
\newcommand{\lambdafull}[0]{\bm{\lambda}}

\DeclareMathOperator{\disp}{\mathbf{u}}

\newcommand{\gravity}[0]{\vec{g}}
\newcommand{\rhoRef}[0]{\rho^f_\text{ref}}
\newcommand{\LameName}{Lam\'e\xspace}
\newcommand{\PecletName}{Pecl\'et\xspace}
\newcommand{\specHeatCapacity}[0]{{c^f}}
\newcommand{\compres}[0]{{\gamma}}
\renewcommand\vec{\mathbf}

\crefname{algorithm}{algorithm}{algorithms}
\Crefname{algorithm}{Algorithm}{Algorithms}
\crefname{ineq}{inequality}{inequalities}
\Crefname{ineq}{Inequality}{Inequalities}

\newcommand{\email}[1]{\href{mailto:#1}{#1}}

\title{A block preconditioner for thermo-poromechanics with frictional deformation of fractures}

\date{}
\hypersetup{
pdftitle={A block preconditioner for thermo-poromechanics with frictional deformation of fractures},
pdfsubject={A block preconditioner for thermo-poromechanics with frictional deformation of fractures},
pdfauthor={Yury Zabegaev, Inga Berre, Eirik Keilegavlen},
pdfkeywords={Thermo-poromechanics, 
Frictional contact mechanics,
Mixed-dimensional,
Preconditioning,
Fixed-stress,
Constrained Pressure Residual (CPR)}
}

\begin{document}

\author{
\href{https://orcid.org/0009-0006-7095-3044}{\includegraphics[scale=0.06]{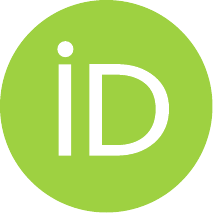}Yury Zabegaev}
\thanks{Center for Modeling of Coupled Subsurface Dynamics, Department of Mathematics, University of Bergen, Bergen, Norway}\\
\email{yury.zabegaev@uib.no}
\and
\href{https://orcid.org/0000-0002-0212-7959}{\includegraphics[scale=0.06]{orcid.pdf}Inga Berre}\footnotemark[1]\\
\email{inga.berre@uib.no}
\and
\href{https://orcid.org/0000-0002-0333-9507}{\includegraphics[scale=0.06]{orcid.pdf}Eirik Keilegavlen}\footnotemark[1]\\
\email{eirik.keilegavlen@uib.no}
}

\maketitle

\begin{abstract}
The numerical modeling of fracture contact thermo-poromechanics is crucial for advancing subsurface engineering applications, including CO\textsubscript{2} sequestration, production of geo-energy resources, energy storage and wastewater disposal operations.
Accurately modeling this problem presents substantial challenges due to the complex physics involved in strongly coupled thermo-poromechanical processes and the frictional contact mechanics of fractures. 
To resolve process couplings in the resulting mathematical model, it is common to apply fully implicit time stepping. This necessitates the use of an iterative linear solver to run the model. 
The solver's efficiency primarily depends on a robust preconditioner, which is particularly challenging to develop because it must handle the mutual couplings between linearized contact mechanics and energy, momentum, and mass balance. In this work, we introduce a preconditioner for the problem based on the nested approximations of Schur complements. To decouple the momentum balance, we utilize the fixed-stress approximation, extended to account for both the porous media and fracture subdomains. The singularity of the contact mechanics submatrix is resolved by a linear transformation. Two variations of the algorithm are proposed to address the coupled mass and energy balance submatrix: either the Constrained Pressure Residual (CPR) or the System-AMG approach. The preconditioner is evaluated through numerical experiments of fluid injection into fractured porous media, which causes thermal contraction and subsequent sliding and opening of fractures. The experiments show that the preconditioner performs robustly for a wide range of simulation regimes governed by various fracture states, friction coefficients and \PecletName number. The grid refinement experiments demonstrate that the preconditioner scales well in terms of GMRES iterations, in both two and three dimensions.
\end{abstract}


\section{Introduction}
Fractures have a significant and manifold impact on multiphysics processes in the subsurface. 
Fractures, and networks of fractures, can form connected pathways that enable fast transport of mass and energy through otherwise low-permeable rocks.
Fractures also deform under a combination of mechanical, hydraulic, and thermal, thermo-poromechanical for short, effects. Such deformation is associated with both alterations of the fracture permeability through dilation or contraction, and with induced seismic events \cite{ellsworth2013injection}.
Subsurface applications for which such dynamics can take place include 
CO\textsubscript{2} sequestration \cite{FAN20191054,LIU20191,flemisch2024fluidflower}, geothermal energy production \cite{PAN201919,WEI2019120,ASAI2019763}, and underground gas storage \cite{zhou2019seismological,karev2019geomechanical,FIRME2019103006}.

Numerical simulations are an essential tool to better understand how fractures deform as in response to subsurface engineering. However, accurately modeling the processes requires  the representation of mutual interactions between  transport of mass and energy in fractures and the host rock, deformation of the host rock due to thermo-poromechanics \cite{coussy2004poromechanics}, and frictional sliding and opening of fractures. 
By themselves, each of these are complex processes, and resolving their combination in a simulation model is highly challenging. 
In particular, Coulomb friction coupled with the non-penetration condition causes non-smooth shifts in constitutive laws as it involves transitions between distinct fracture states: opening, sticking and sliding \cite{kikuchi1988contact,matrins1987,sofonea2012mathematical}. 
As a further complication, the geometry of fracture networks can be complex, with the orientation of individual fractures relative to the background stress field having a first-order impact on the dynamics.

To address these challenges, we adopt the explicit fracture representation in terms of a discrete fracture-matrix modeling approach 
\cite{GARIPOV2019104075,WANG2021103985,WANG2020114693,Brenner2016,bonaldi2024numerical,BONALDI2022110984,stefansson_fully_2021,berre2019flow,porepy2024,berge_finite_2020,porepy2021}, 
where fractures and their intersections are represented explicitly as lower-dimensional geometrical objects embedded within the computational domain.
The explicit representation of fractures allows us to directly model flow and transport inside fractures and interaction between fractures and the host rock.
The strong coupling of multiple physics calls for a fully implicit time scheme, along with the application of Newton's method to solve the resulting discretized system.
Consequently, the complexity of the problem is transferred to the linear solver, which requires an efficient an robust preconditioner.

In our previous work \cite{zabegaev_efficient_2025}, we developed a preconditioner for problems that couple poromechanics in fractured porous media with fracture deformation governed by contact mechanics.
Our approach was based on the so-called fixed-stress approximation \cite{white_block-partitioned_2016,mikelic_convergence_2013,kim_stability_2011,storvik_optimization_2018,BOTH2017101}
to decoupled the poromechanics problem, following the extension to fractured porous media suggested in \cite{girault_convergence_2016}.
To overcome the potential non-degeneracy of the contact mechanics submatrix was invertible, we applied a linear transformation originally suggested in saddle-point theory \cite{bacq_allatonce_2023,notay_new_2016,benzi_numerical_2005}.

New challenges arise for the non-isothermal extension of the problem: 
By itself, the energy conservation equation includes an elliptic and a hyperbolic component that both must be accounted for in a preconditioning strategy. 
Furthermore, the energy conservation equation will be mutually coupled with the conservation of both fluid mass and momentum.
Of these couplings, the former can be expected to be by far the strongest under normal subsurface conditions, and it is also much studied in the context of classical reservoir simulation.

Most non-isothermal models in reservoir simulation employ the two-stage Constrained Pressure Residual (CPR) preconditioner \cite{wallis_incomplete_1983,wallis_constrained_1985} to treat the coupling between fluid flow and transport. 
Originally developed for multi-component isothermal models, CPR effectively addresses the elliptic structure of the flow equation using AMG and handles the hyperbolic transport effects by applying an incomplete LU factorization (ILU) to the coupled transport and flow problem. 
While CPR is a robust solution for multi-component isothermal transport, its performance deteriorates if
transport is dominated by diffusion which ILU is not well suited to approximate. 
As shown in \cite{roy_block_2019, cremon_constrained_2024}, such diffusion dominated regimes can be relevant for problems involving thermal transport.
A different approach addresses this issue by applying the System-AMG method \cite{boomeramg,gries2015system}, which applies the AMG method to the monolithic pressure-temperature submatrix \cite{cremon_constrained_2024}. Alternatively, the flow and energy equations can be decoupled using one more Schur complement approximation \cite{roy_block_2019,cremon_multi-stage_2020}, or by reinforcing the CPR method with one more AMG application for the energy equation \cite{mohajeri_novel_2020}.
For our target simulations, the energy transport both in fractures and in the surrounding matrix amplifies the difference between dominant thermal effects in different (spatial) regions of the model. 
It is therefore not clear \textit{a priori} whether a standard CPR-based approach is sufficiently effective, or if more elaborate preconditioners are needed for the flow-energy subproblem.

The solution strategies for the coupled energy balance and poromechanics equations are less extensively covered in the literature. Various extensions of the fixed-stress approximation have been suggested, resulting in development of the convergent sequential iterative schemes for the thermo-poromechanics problem, with the potential to recast them as preconditioners for the fully coupled problem \cite{both_gradient_2019,kim_unconditionally_2018,brun_monolithic_2020}. 
Significant progress has been made in the development of preconditioners for coupled poromechanics and mass transfer \cite{white_two-stage_2019,bui_scalable_2020}, which share a similar structure in the to our problem.
The elasticity field is decoupled first by using the extended fixed-stress approximation to build the Schur complement of the fluid flow and transport subproblem, which is treated by CPR. The similar approach is developed in the recent work \cite{novikov_finite_2024}, where the heat and mass transport are treated monolithically within the CPR framework, while the fixed-stress approximation is used to decouple the mechanics. These results suggest the appropriate decoupling order in the preconditioner for our problem. 
However, to the best of the authors' knowledge, no CPR-like preconditioner for the flow and transport problem both in the porous medium and in explicitly represented fractures, including the transport between them, has been reported. 
Likewise, preconditioners that also account for the coupling with fracture contact mechanics are also not available.

In this paper, we propose two variations of the preconditioner that effectively combine approaches to address the contact mechanics, poromechanics and energy transport in a fractured porous medium. We apply a linear transformation to eliminate the singularity caused by the saddle-point formulation of the contact mechanics and construct a series of nested Schur complement approximations to decouple the full problem into a sequence of manageable subsystems. The elasticity equation is decoupled using the fixed-stress approximation, and
two methods for solving the coupled pressure-temperature subproblem are explored: either CPR or System-AMG. We probe the robustness of the preconditioner with these variations in numerical experiments.

The paper is structured as follows: \Cref{sec:model_summary} gives a brief overview of the mathematical model, while the full description including governing equations, constitutive laws and boundary conditions is given in \Cref{appendix:model_description}. \Cref{sec:discretization} covers the discretization of the model.
\Cref{sec:preconditioner} describes the preconditioner, consisting of the linear transformation to handle the saddle point structure, the fixed-stress approximation to decouple the elasticity, and the CPR and the System-AMG alternatives to solve the pressure-temperature subproblem. In \Cref{sec:numerical_experiments}, we present numerical experiments to evaluate the performance of the preconditioner. The concluding remarks are provided in \Cref{sec:results}.

\section{Model overview}
\label{sec:model_summary}

This section provides a brief summary of the mathematical model for mixed-dimensional flow, deformation, and energy transport in a thermo-poroelastic medium with frictional fracture deformation, governed by contact mechanics. 
The model was introduced previously in \cite{stefansson_fully_2021,porepy2024}.
Since our focus herein is not on physical modeling, we provide 
a brief summary of the model, emphasizing the features that are important for deriving the preconditioners.
The full model is listed in \Cref{appendix:model_description}, and we refer to equations there when relevant.

\begin{wrapfigure}{R}{7.5cm}
    \centering
    \includegraphics[width=7.5cm]{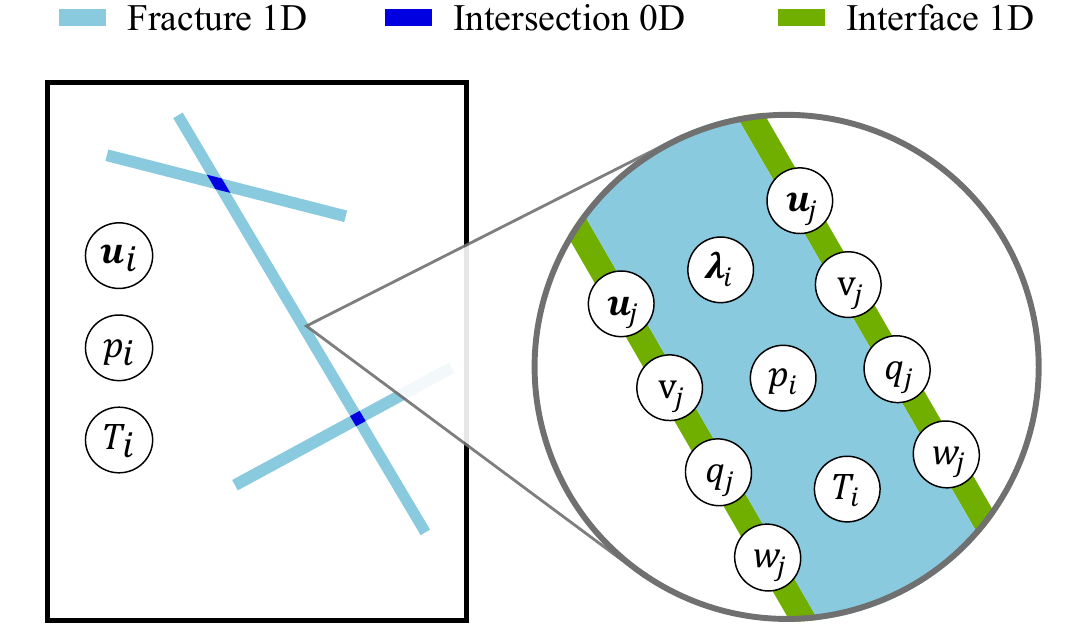}
    \caption{Rectangular 2D porous medium with three embedded 1D fractures and two 0D intersections.
    Symbols inside the circles represent the primary variables and their domains of definition. 
    The magnifying glass on the right shows the variables defined within the fracture and on the interfaces between the fracture and the porous medium.}
    \label{fig:primary_variables}
\end{wrapfigure}

The model geometry consists of a 2D or 3D deformable porous medium with embedded fractures. These fractures are respectively 1D or 2D geometric objects, usually referred to as co-dimension 1 objects. Fractures can intersect to form objects of co-dimension 2 and co-dimension 3 (only in 3D). 
Together, the porous medium and the individual lower-dimensional objects are collectively referred to as subdomains. 
The interaction between the lower- and higher-dimensional objects is facilitated through interfaces, which are the geometrical objects representing the fracture boundaries. The example model geometry is illustrated in \Cref{fig:primary_variables}.

Below we summarize the governing equations,
see also \Cref{appendix:model_description} for the detailed model:
\begin{enumerate}
    \item Deformation of fractures is described by inequalities that impose a non-penetration condition in the normal direction \eqref{ineq:contact_normal_ineq} and a Coulomb friction law \eqref{ineq:contact_tangential_ineq} in the tangential direction. 
    In practice, the inequalities are recast to equalities through the augmented Lagrangian approach and \cref{eq:contact_equalities} will be present in our system of equations; 
    \item Balance of forces across is imposed on the interfaces between fractures and the porous medium, according to \cref{eq:intf_force_balance};
    \item Quasi-static momentum balance, \cref{eq:momentum_conservation} is imposed in the porous medium;
    \item Inter-dimensional fluxes caused by energy diffusion, energy convection and fluid mass are governed by \cref{eq:interface_flow,eq:interface_energy_advective,eq:interface_energy_diffusive}, respectively. These equations are imposed both on the interfaces between fractures and matrix, and between fractures and fracture intersections;
    \item The fluid mass balance \eqref{eq:mass_conservation} is defined in all subdomains;
    \item The energy balance \eqref{eq:energy_conservation} is defined in all subdomains.
\end{enumerate}

We also list the primary variables of our model, see also \Cref{fig:primary_variables}:
\begin{enumerate}
    \item The fracture traction, $\lambdafull_i$, is defined in the fracture subdomains;
    \item The fracture interface displacement, $\disp_j$, is defined on the interfaces between fractures and the porous medium;
    \item The displacement, $\disp_i$, is defined in the porous medium;
    \item The inter-dimensional fluxes associated with energy diffusion, $q_j$, energy convection, $w_j$, and fluid mass flux $v_j$, are all defined on all interfaces;
    \item The fluid pressure, $p_i$, is defined in all subdomains;
    \item The temperature, $T_i$, is defined on all subdomains.
\end{enumerate}
The subscripts distinguish the unknowns defined on subdomains ($i$) from those on interfaces ($j$). 
The ordering of equations and variables in the above lists will be reflected in the linear system and preconditioner studied in \Cref{sec:preconditioner}.
Dirichlet or Neumann boundary conditions are applied on the external boundaries, while continuity of displacements and normal fluid and energy fluxes are enforced as boundary conditions at the interfaces.
According to the formulation, three local fracture states are allowed: its sides can either stick, slide, or open relative to each other.

The constitutive laws for the thermo-poroelastic medium are in large part based on Coussy \cite{coussy2004poromechanics}, here we list the main laws. 
The fluid is considered slightly compressible with the fluid flux governed by a Darcy's law.
The permeability is taken as constant in the porous medium, while in the fractures, the permeability depends on the aperture through the cubic law \cite{zimmerman1996hydraulic}. Differences in permeability between fractures and porous medium can lead to strong local variations in flux magnitudes.
The energy conservation equation includes both the diffusive and convective transport terms.
The thermo-poroelastic coupling is established through multiple dependencies: stress is a function of displacement, pressure, and temperature; fluid density depends on pressure and temperature; and porosity depends on pressure, temperature, and displacement according to \cref{eq:poromech_stress,eq:fluid_density,eq:porosity} in \Cref{appendix:model_description}.

In terms of mathematical structure, which is critical for the construction of preconditioners,
the fluid mass balance equation is parabolic, while the momentum balance is elliptic.
The presence of convection and diffusion implies that the energy conservation equation has a mixed parabolic-hyperbolic character; potentially with different physics dominating in different parts of the domain.
Finally, we remark that the augmented Lagrangian formulation of contact mechanics introduces a saddle-point structure into the system.

\section{Discretization}
\label{sec:discretization}

This section summarizes the main ingredients of the discretization of the governing equations, %
see \cite{stefansson_fully_2021,porepy2024,porepy2021} for the details.
The domain is discretized using a simplex grid that aligns with fractures and intersections, with grid construction handled by Gmsh \cite{gmsh}. This approach ensures that each lower-dimensional subdomain geometrically coincides with a set of faces in the surrounding higher-dimensional grid.
The spatial discretization is based on a family of cell-centered finite volume schemes. The multi-point flux approximation (MPFA) method \cite{Aavatsmark2002} is applied to discretize the diffusive scalar fluxes (mass and energy), whereas the multi-point stress approximation (MPSA) method \cite{Nordbotten2016} is used to discretize the thermo-poromechanical stress.
Coupling of discretizations on different subdomains follows the approached defined in \cite{nordbotten2019unified}.
The first-order upwinding scheme is utilized to discretize the advective fluxes of energy and inter-dimensional fluid flux; the values of the fluid flux $\vec{v}_i$ and the inter-dimensional fluid flux $v_j$ are computed from the solution at the previous Newton iteration. The temporal discretization is performed using the implicit Euler scheme.

The discretized governing equations from \Cref{sec:model_summary} result in a nonlinear system that must be solved to advance the simulation in time. Due to the contact mechanics relations, the system is only semi-smooth.
To solve it, we employ a semi-smooth Newton method \cite{ito2003semi}, which involves inverting the generalized Jacobian matrix. This matrix coincides with the original Jacobian in smooth regions and adopts derivative values from one side or the other in non-smooth regions \cite{ito2003semi,berge_finite_2020}. Hereafter, we refer to it simply as the Jacobian, $J$.

\section{Preconditioner variants for mixed-dimensional contact thermo-poromechanics}

\label{sec:preconditioner}
In this section, we construct two variations of the preconditioning algorithm for the Jacobian of the linearized problem
\begin{equation}
    J\vec{x}=\vec{r}\,,
\end{equation}
where $x$ is the unknown vector and $r$ is the residual right hand side.
The preconditioner in its core is based on decoupling the Jacobian $J$ into a sequence of single-physics subproblems that can be handled by tailored preconditioners.
The main procedure is the block $\mathcal{LDU}$ factorization of a $2\times2$ block matrix:
\begin{equation}
\label{eq:example_block_mat_2x2}
J=
    \begin{bmatrix}
        A & B \\ C & D
    \end{bmatrix}
    = \mathcal{L D U} =
    \begin{bmatrix}
        I              &  \\
        C  A^{-1} & I
    \end{bmatrix}
    \begin{bmatrix}
        A & \\
          & S_A
    \end{bmatrix}
    \begin{bmatrix}
        I & A^{-1} B \\
          & I
    \end{bmatrix},
\end{equation}
where $S_A$ is the Schur complement: $S_A \coloneq D - C A^{-1} B$.
The inverted upper-triangular matrix $(\mathcal{DU})^{-1}$ can be used as a right preconditioner, leading to
$\mathcal{LDU}(\mathcal{DU})^{-1}=\mathcal{L}$. Since the diagonal elements of $\mathcal{L}$ all equal one, so do its eigenvalues, thus a Krylov subspace method such as GMRES combined with this preconditioner should converge in two iterations assuming exact arithmetics. In practice, the exact inverse $(\mathcal{DU})^{-1}$ cannot be constructed within a reasonable computational cost, and instead the approximate operator is computed:
\begin{equation}
    (\mathcal{DU})^{-1} \approx \begin{bmatrix}
        I & -\tilde{A}^{-1} B \\
          & I
    \end{bmatrix}
    \begin{bmatrix}
        \tilde{A}^{-1} &    \\
             & \tilde{S}_A^{-1}
    \end{bmatrix},
\end{equation}
where $\tilde{A}^{-1}$ and $\tilde{S}_A^{-1}$ are the approximations of the corresponding inverted matrices. 
Moreover, as the Schur complement $S_A$ cannot be explicitly constructed, the approximation $\tilde{S}_A^{-1}$ must be based on a sparse approximation of the Schur complement.
What are suitable inexact subsolvers $\tilde{A}^{-1}$ and $\tilde{S}_A^{-1}$ and Schur complement approximation $\tilde{S}_A$ depend on the underlying physics. 
For block matrices with more than two block rows and columns, a recursive $\mathcal{LDU}$ decomposition applies, and the corresponding block-triangular preconditioner is a suitable option. In this case, nested Schur complements must be approximated together with their inverses.

\begin{figure}[htbp]
    \centering
    \includegraphics[width=1\linewidth]{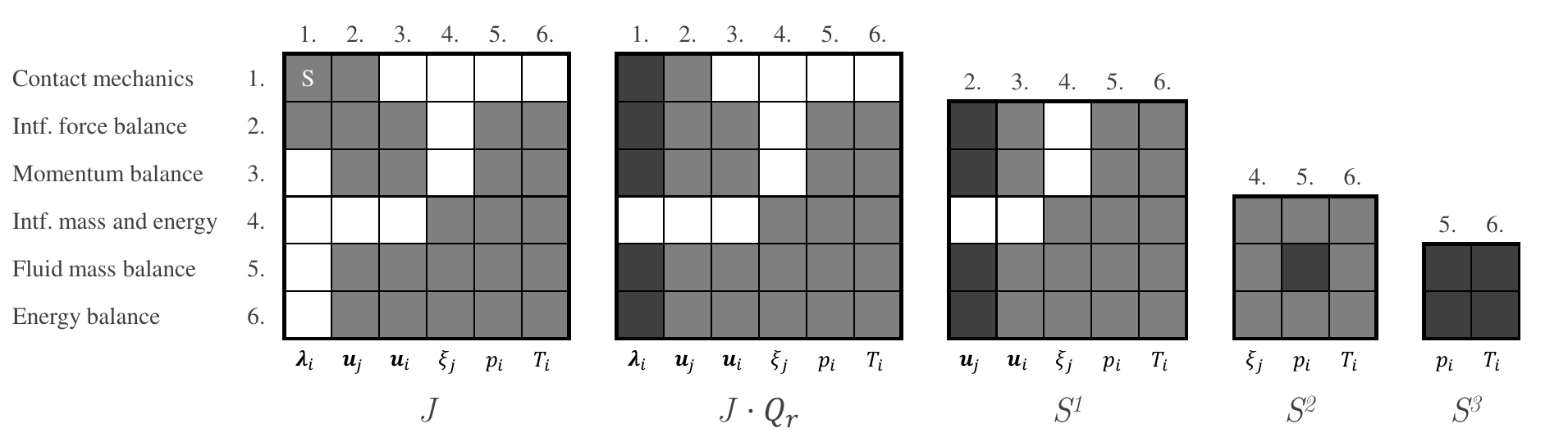}
    \caption{The block structure of the matrices at different stages of the preconditioner algorithm. White cells indicate empty submatrices, light-gray cells represent non-empty submatrices from the original Jacobian, and dark-gray cells denote modified submatrices. The letter ``S'' marks a singular submatrix. From left to right, the stages include: the original Jacobian, the Jacobian after the linear transformation (\Cref{sec:linear_transformation}), the first-level Schur complement (\Cref{sec:eliminating_contact_mechanics}), the second-level Schur complement (\Cref{sec:eliminating_force_balance}), and the third-level Schur complement (\Cref{sec:eliminating_intf_flow}). Column 4 of the matrix refers to all the interface fluxes: $\xi_j \in \{ v_j, w_j, q_j\}$.
    }
\label{fig:elimination_order}
\end{figure}

We schematically represent the Jacobian of our problem in \Cref{fig:elimination_order} (left). The rows and the columns of the block matrix correspond to the specific governing equations and the unknowns introduced in \Cref{sec:model_summary}, following the same order. Row 4 combines the equations of the interface fluid and energy fluxes \cref{eq:interface_flow,eq:interface_energy_advective,eq:interface_energy_diffusive}, while column 4 refers to all the interface flux unknowns: $\xi_j \in \{ v_j, w_j, q_j\}$.
We will refer to the submatrices of the Jacobian as $J_{ij}$ with subscripts $i, j \in \{ 1, 2, 3, 4, 5, 6 \}$ according to the illustration, e.g. $J_{11}$ corresponds to the contact mechanics submatrix.
For our problem, the contact mechanics submatrix $J_{11}$ is singular and corresponds to the saddle-point structure of the matrix. It is strongly coupled with the elasticity submatrix, which, in turn, is strongly coupled with the fluid and the energy mass balance submatrices.

As in our previous work on the isothermal problem \cite{zabegaev_efficient_2025}, we use a fixed-stress based approximation to decouple elasticity from fluid flow, along with a linear transformation to handle the saddle-point structure caused by the contact mechanics submatrix. 
The key novelty of this article lies in the treatment of the energy equation, 
but we nevertheless provide an overview of all preconditioner steps below.

The algorithm consists of a preprocessing stage that handles singularities in the contact mechanics submatrix as detailed in \Cref{sec:linear_transformation}, followed by an application stage. 
This second stage involves the construction of a sequence of linear subproblems through a decoupling of the original problem, and of solving the resulting problems.
For $S^i$ denoting the $i$-th Schur complement, the steps are the following:
\begin{enumerate}
\item Decouple the contact mechanics equations and recast the remainder as $S^1$ as described in \Cref{sec:eliminating_contact_mechanics};
\item Decouple the momentum balance equation from $S^1$ and recast the remainder as $S^2$ through the fixed-stress based approximation, see \Cref{sec:eliminating_force_balance};
\item Decouple the interface fluid and energy fluxes from $S^2$ and recast the remainder as $S^3$, see  \Cref{sec:eliminating_intf_flow};
\item Approximately inverse the coupled pressure-temperature submatrix $S^3$ using one of the variations proposed in \Cref{sec:eliminating_energy}.
\end{enumerate}
The order of decoupling of equations and the sparsity pattern of the reduced matrices is illustrated in \Cref{fig:elimination_order}.

\subsection{Preprocessing}
\label{sec:linear_transformation}
The contact mechanics submatrix $J_{11}$ consists of block-diagonal square matrices of size $2$ or $3$, depending on the spatial dimension, with one block for each cell.
Blocks that correspond to fracture cells in a sliding or sticking cells may be singular.
To facilitate the elimination of $J_{11}$, 
we observe that the contact mechanics equations are coupled only with the interface force balance equation, and the Jacobian matrix can be expressed as
\begin{equation}
\label{eq:jac_contact}
    J = \begin{bmatrix}
        J_{11} & J_{12} & 0 \\
        J_{21} & J_{22} & J_{2 \star} \\
        0               & J_{\star 2} & J_{\star \star} \\
    \end{bmatrix},
\end{equation}
where the symbol $\star$ represents the remainder, specifically the submatrices with indices $i, j \in \{3, 4, 5, 6\}$.
This serves as motivation for a linear transformation matrix, inspired by \cite{notay_new_2016,bacq_allatonce_2023}:
\begin{equation}
    \label{eq:Qright}
    Q_r = \begin{bmatrix}
        I_{11}& 0& 0 \\
        -D_{22}^{-1}  J_{21}& I_{22}& 0 \\
        0& 0& I_{\star \star} \\
    \end{bmatrix},
\end{equation}
where the subscript $r$ marks that $Q_r$ is applied from the right and $I_{ii}$ is an identity matrix with corresponding to the block structure of the Jacobian.
${D_{22}}$ represents the block diagonal of $J_{22}$, with block size $2$ or $3$ depending on the spatial dimension, hence $D_{22}^{-1}$ can be computed at a low computational cost.
Defining $\tilde{J} \coloneq J  Q_r$ and $\tilde{\vec{x}} \coloneq Q_r^{-1} \vec{x}$, we can solve the transformed linear system $\tilde{J}  \tilde{\vec{x}} = \vec{r}$. 
The transformed matrix, shown second from the left in \Cref{fig:elimination_order}, is assembled explicitly. 
It is important to note that the transformation is applied only once, during preprocessing.
The sparsity structure of the Jacobian block matrix after the linear transformation becomes different:
\begin{equation}
\label{eq:J_tilde}
    \tilde{J} = \left[\begin{array}{c;{4pt/2pt}cc}
        J_{11} - J_{12}  D_{22}^{-1}  J_{21} & J_{12} & 0 \\   \hdashline[4pt/2pt]
        E_{21} & J_{22} & J_{2 \star} \\
        -J_{\star 2}  D_{22}^{-1}  J_{21} & J_{\star 2} & J_{\star \star}
    \end{array}\right],
\end{equation}
where $E_{21} \coloneq \left( I_{22} - J_{22}  D_{22}^{-1} \right) J_{21}$ is nonzero due to the inexact inverse of $J_{22}$.

\subsection{Contact mechanics subsolver and the first-level Schur complement}
\label{sec:eliminating_contact_mechanics}
After the linear transformation, the submatrix $\tilde{J}_{11} \coloneq J_{11} - J_{12}  D_{22}^{-1} J_{21}$ retains its block-diagonal structure of square matrices of size 2 or 3 (spacial dimension), but now all of them are invertible. Inversion can be done inexpensively to produce the exact first-level Schur complement
\begin{equation}
\label{eq:S^1}
    S^1 = \left[\begin{array}{cc;{4pt/2pt}ccc}
        S^1_{22} & J_{23} &        & J_{25} & J_{26} \\
        S^1_{32} & J_{33} &        & J_{35} & J_{36} \\ \hdashline[4pt/2pt]
                 &        & J_{44} & J_{45} & J_{46} \\
        S^1_{52} & J_{53} & J_{54} & J_{55} & J_{56} \\
        S^1_{62} & J_{63} & J_{64} & J_{65} & J_{66} \\
    \end{array}\right],
\end{equation}
where $S^1_{22} = J_{22} + E_{22}$, $S^1_{i2} = J_{i2} - \tilde{J}_{i1} \cdot \tilde{J}_{11}^{-1} \cdot J_{12}$ for $i \in \{ 3, 4, 5, 6\}$, and $E_{22} = -E_{21} \tilde{J}_{11}^{-1} J_{12}$.
The Schur complement $S^1$ is depicted as matrix 3 in \Cref{fig:elimination_order}.
Note that the submatrices in all columns except column 2 remain unchanged.

\subsection{Momentum balance subsolver and the second-level Schur complement}
\label{sec:eliminating_force_balance}
The second Schur complement involves the elimination of the momentum balance equations.
The submatrix $\{2,3\}$, which represents the coupled interface force balance and momentum balance operator, is treated monolithically by applying AMG. We decouple it from the remainder of the system by constructing the second-level Schur complement $S^2$.
As can be seen from \eqref{eq:S^1}, there is no coupling between mechanics (blocks 2 and 3) and the interface fluid/energy fluxes (block 4). However, the Schur complement for submatrices 5 and 6 still needs to be approximated:
\begin{equation}
\label{eq:S^2}
    S^2 = \left[\begin{array}{c;{4pt/2pt}cc}
         J_{44} & J_{45} & J_{46} \\
         \hdashline[4pt/2pt]
         J_{54} & J_{55} & J_{56} \\
         J_{64} & J_{65} & J_{66}
    \end{array}\right] 
    +
    \left[\begin{array}{c;{4pt/2pt}cc}
    0 & 0 & 0 \\
    \hdashline[4pt/2pt]
    0 & \mathcal{F}_{55} & \mathcal{F}_{56} \\
    0 & \mathcal{F}_{65} & \mathcal{F}_{66}
    \end{array}\right]. 
\end{equation}
The submatrices $\mathcal{F}_{ij}$ correspond to the second term in the Schur complement formula and represent the decoupling of the momentum balance equation from the pressure-temperature subproblem. 
For the poromechanics coupling, we apply the fixed-stress based approximation as follows:
\begin{equation}
\label{eq:fixed_stress_S_55}
    \mathcal{F}_{55} \approx \frac{\rho_i V_i}{\Delta t} L_i,
\end{equation}
where $L_i$ is the fixed-stress stabilization matrix, divided by the simulation time step $\Delta t$ and multiplied by the fluid density $\rho_i$ times the (mixed-dimensional) cell volume $V_i$. The subscript $i$ corresponds to a specific subdomain, as the expression for the fixed-stress stabilization matrix varies depends on the dimensionality of the subdomains \cite{mikelic_convergence_2013, kim_stability_2011, white_block-partitioned_2016, girault_convergence_2016, zabegaev_efficient_2025}. For the porous medium dimension and the fractures, the fixed-stress stabilization matrices $L_\text{pm}$ and $L_\text{frac}$ are defined as follows:
\begin{equation}
    \label{eq:fixed_stress_coefs}
    L_\text{pm} = \dfrac{\alpha^2}{2 G / D + \Lame};
    \quad
    L_\text{frac} = \dfrac{\ujumpn \alpha^2 \compres}{\Lame (1/M + \varphi \compres)},
\end{equation}
where $\dfrac{1}{M} = \dfrac{\partial \varphi}{\partial p} = \dfrac{(\alpha - \varphi_\text{ref})(1 - \alpha)}{\Lame + \frac{2}{3}G}$, $\alpha$ is the Biot coefficient, $G$ and $\Lame$ are the Lame coefficients, $D$ is the ambient dimension of the problem, $\varphi$ is the porosity, $\gamma$ is the fluid compressibility and $\ujumpn$ is the normal displacement jump between fracture sides. The first coefficient refers to the classical fixed-stress formulation \cite{mikelic_convergence_2013, kim_stability_2011, white_block-partitioned_2016}, while the second is derived in \cite{girault_convergence_2016} for the fractured poromechanics problem without contact mechanics. It was then tested for the isothermal version of our problem in \cite{zabegaev_efficient_2025}.
Since the intersection subdomains do not interact with the momentum balance equation, they remain uncoupled and, therefore, do not require any stabilization.

Stabilization terms to decouple the fracture-less mechanics-energy system are studied in \cite{both_gradient_2019,kim_unconditionally_2018,brun_monolithic_2020}, and the following approximation is used:
\begin{equation}
\mathcal{F}_{66} \approx \frac{\rho V}{\Delta t} L_\text{thermal};
\quad
L_\text{thermal} = {c} \frac{(\beta^s)^2}{2 G / D + \Lame},
\end{equation}
where ${c}$ is a unit-less constant and $\beta^s$ is the volumetric solid thermal expansion. To the authors' knowledge, no similar stabilization coefficient for the energy transport in fractures has been developed.
For the material parameters of interest, we observe that the values in the matrix resulting from these stabilization terms are several orders of magnitude smaller than $J_{66}$, especially when compared to the difference in magnitudes between the terms of the poromechanical fixed-stress stabilization and $J_{55}$.
Numerical experiments indicate negligible difference in linear solver iterations with or without thermal stabilization, as well as in exactly inverted Schur complements for small problems (results not shown here). Consequently, we conclude that the mechanics-energy coupling is significantly weaker than the poromechanical coupling and does not require specific stabilization for the type of problem we are considering. Hence, we set:
\begin{equation}
    \mathcal{F}_{56} \approx 0;
    \quad
    \mathcal{F}_{66} \approx 0;
    \quad
    \mathcal{F}_{65} \approx 0.
\end{equation}

\subsection{Third-level Schur complement: Eliminating the interface fluxes}
\label{sec:eliminating_intf_flow}

Row 4 of the block matrix representation of the Jacobian correspond to the interface fluxes of mass and energy. The approach to approximate the inverse of the interface fluid mass flux submatrix with the diagonal was applied in \cite{hu_effective_2023,zabegaev_efficient_2025}. We extend this idea and utilize the diagonal approximation to invert the matrix $J_{44}$:
\begin{equation}
\label{eq:S^3}
    S^3 = 
    \begin{bmatrix}
        S^3_{55} & S^3_{56} \\
        S^3_{65} & S^3_{66} \\
    \end{bmatrix}
    =
    \begin{bmatrix}
        S^2_{55} & J_{56} \\
        J_{65} & J_{66} \\
    \end{bmatrix}
    -
    \begin{bmatrix}
        J_{54} \\ J_{64}
    \end{bmatrix}
    \cdot 
    \text{diag}(J_{44})^{-1} 
    \cdot
    \begin{bmatrix}
        J_{45} & J_{46}
    \end{bmatrix},
\end{equation}
where $\text{diag} (\cdot)^{-1}$ is an operator that constructs a diagonal matrix based on the inverted diagonal of its input matrix.
The linear system based on $J_{44}$ is solved by an ILU(0) algorithm.

\subsection{Energy and mass conservation subsolvers}
\label{sec:eliminating_energy}

We define the preconditioner for the coupled energy and mass conservation equation as the operator $M^{-1}$, which approximates the inverse of $S^3$. Due to the mixed parabolic-hyperbolic nature of the energy equations, two options are considered.

\paragraph{Constrained Pressure Residual (CPR)}
The standard approach to decouple the hyperbolic transport from the pressure field is the CPR preconditioner \cite{wallis_incomplete_1983,wallis_constrained_1985,roy_block_2019,roy_constrained_2020,cremon_multi-stage_2020,cremon_constrained_2024}. This is a two-stage preconditioner, which consists of the local and the global preconditioners $M_G^{-1}$ and $M_L^{-1}$, respectively. The action of a two-stage preconditioner for a generic matrix $A$ is given by:
\begin{equation}
\label{eq:two_stage_prec}
    M^{-1} = M_G^{-1} + M_L^{-1} (I - A M_G^{-1}).
\end{equation}
The global preconditioner targets the long-distance interaction in the elliptic pressure field. Using the same the block matrix structure as in $S^3$, the global preconditioner is:
\begin{equation}
    M_G^{-1} = \begin{bmatrix}
        M_{55}^{-1} & 0 \\
        0         & 0 
    \end{bmatrix},
\end{equation}
where $M_{55}^{-1}$ denotes the AMG application to $S^3_{55}$.
The local preconditioner applies ILU(0) to the coupled pressure-temperature submatrix $S^3$. For ILU(0) applications, the matrix is rearranged to group the variables and the equations by cells rather than by physics (interleaved). We note that the original CPR algorithm for the multi-component mass transfer problem includes a decoupling step in the global preconditioner, but it was shown practical to skip it for the energy transfer problem \cite{roy_block_2019,roy_constrained_2020}. The CPR algorithm effectively tackles the hyperbolic energy transfer, but can struggle when diffusion becomes predominant in the energy equation \cite{roy_block_2019,cremon_multi-stage_2020}.

\paragraph{System-AMG}
The system-AMG approach applies AMG to the coupled pressure-temperature matrix $S^3$ \cite{roy_constrained_2020,cremon_constrained_2024} and can serve a good alternative for CPR. The efficiency of this approach depends on how well the underlying AMG implementation can handle the coupled system. Notably, AMG works best for diffusion-driven problems, while hyperbolic effects may deteriorate its performance.
It is not clear \textit{a priori} whether this will have significant impact on for practical problems.

\subsection{Preconditioner algorithms summary}
The proposed variations of the preconditioner approximate the inverse of the upper-triangular block matrix of the block $\mathcal{LDU}$ factorization: $P^{-1} \approx (\mathcal{DU})^{-1}$, where the block factorization is applied to the linearly transformed Jacobian. They can be generalized by the following block matrix:
\begin{equation}
P = \left[\begin{array}{cccccc}%
\tilde{J}_{11} & J_{12}  &       &       &       &       \\ %
               & S^1_{22}& J_{23}&       & J_{25}& J_{26}\\
               & S^1_{32}& J_{33}&       & J_{35}& J_{36}\\ %
               &         &       & J_{44}& J_{45}& J_{46}\\ %
               &         &       &       & M_{55}& M_{56}\\
               &         &       &       & M_{65}& M_{66}
\end{array}\right],
\end{equation}
where the operator $M_{ij}$ for $i,j \in \{5,6\}$ approximates $S^3_{ij}$ with one of the variations in \Cref{sec:eliminating_energy}. The linearly transformed Jacobian is computed explicitly, and the following linear system is solved by GMRES: $J \cdot Q_r \cdot P^{-1} \cdot \tilde{\vec{x}} = \vec{r}$ for a given right-hand side vector $\vec{r}$. The original vector of unknowns is recovered through an inverse transformation: $\vec{x} = Q_r \cdot P^{-1} \cdot \tilde{\vec{x}}$. The approximate multiplication of $P^{-1}$ by a vector is computed according to \Cref{algorithm:application}. AMG and ILU0 in the algorithm refer to the application of their respective methods. The inverse of the matrix $M_{ij}$ is approximated according to one of the methods in \Cref{sec:eliminating_energy}.
In theory, with exact arithmetic, 
the GMRES method equipped with this preconditioner would converge after just two iterations. However, in real-world applications, its performance is influenced by the accuracy of the approximations of the AMG and ILU subsolvers, as well as the quality of the nested Schur complement approximations. These factors are analyzed in the numerical experiments presented in the next section.

\begin{algorithm}[htbp]
\caption{Preconditioner $P$ application.}
\textbf{Input:} $\vec{r} = [\vec{r}_1, \vec{r}_2, \vec{r}_3, \vec{r}_4, \vec{r}_5, \vec{r}_6]$;
\label{algorithm:application}
\begin{algorithmic}[1]
    \State $\vec{v_5} = \begin{bmatrix}
        M_{55}& M_{56} \\
        M_{65}& M_{66}
    \end{bmatrix}^{-1} \cdot \begin{bmatrix}
         \vec{r_5}  \\
         \vec{r_6}
    \end{bmatrix}$;
    \Comment{Fluid mass and energy balance}
    \State $\vec{v_4} = \text{ILU0}(J_{44}, \vec{r_4} - J_{45} \cdot \vec{v_5} - J_{46} \cdot \vec{v_6})$;
    \Comment{Interface fluid mass and energy flow}
    \State $
    \begin{bmatrix}
        \vec{v_2} \\ \vec{v_3} 
    \end{bmatrix} = \text{AMG}\left(
    \begin{bmatrix}
        S_{22}^1 & J_{23} \\ S_{32}^1 & J_{33}
    \end{bmatrix},
    \begin{bmatrix}
        \vec{r}_2 - J_{25} \cdot \vec{v}_5 - J_{26} \cdot \vec{v_6} \\ \vec{r}_3 - J_{35} \cdot \vec{v}_5 - J_{36} \cdot \vec{v_6}
    \end{bmatrix} \right)$;
    \Comment{Momentum and interface force balance}
    \State $\vec{v}_1 = \tilde{J}_{11}^{-1} \cdot (\vec{r}_1 - J_{12} \cdot \vec{v}_2)$; \Comment{Contact mechanics}
\end{algorithmic}
\textbf{Output:} $[\vec{v}_1, \vec{v}_2, \vec{v}_3, \vec{v}_4, \vec{v}_5, \vec{v}_6]$.
\end{algorithm}

\section{Numerical experiments}
\label{sec:numerical_experiments}
To show the performance of the preconditioners, we present numerical experiments in 2D and 3D domains with fracture networks.
For all experiments, the initial conditions are taken as steady state, and we then consider injection of cold water through an injection well which is located in a fracture; the exact configuration is described in the respective 2D and 3D setups.
The injection will trigger deformation through a combination of two effects: 
On a short time scale, the elevated pressure due to injection will cause a front of higher pressure that propagates rapidly through the fracture network, potentially causing sliding of critically stressed fractures, as well as fracture opening.
On a longer time scale, thermal contraction of the host rock due to cooling will reduce the effective normal forces or fracture surfaces, again causing sliding and potentially opening.

The material constants used in the simulation are listed in \Cref{tab:material_params}.
The reference pressure and temperature are imposed as Dirichlet boundary conditions on all the external boundaries. 
During the simulation, cold fluid is injected with a temperature of $313$ K.
The domain is compressed by applying normal forces on all the external boundaries except one, while the latter remains fixed. The values of the compressive forces are specified for each experiment.
We ignore the effect of gravity.

To improve the conditioning of the linear system, we have found it useful to apply the following scaling of variables and equations:
The unit of mass is set to $10^{10}$ kg, the unit of the contact traction variable $\lambdafull_i$ is scaled by Young modulus, and the opening fracture state indicator is scaled adaptively following the approach outlined in \cite{stefansson2024linesear}.
The energy balance equation in lower-dimensional subdomains is scaled by the subdomain specific volume, defined in \Cref{appendix:model_description}.

An adaptive time-stepping scheme is applied, trying to keep the number of Newton iterations between 3 and 8. The initial time step is on the order of milliseconds and gradually increasing it to the scale of days.
The Newton solution is accepted with the residual relative tolerance of $10^{-7}$.
The line search algorithm is used at Newton iterations.

GMRES is used as the linear solver with the preconditioners described in \Cref{sec:preconditioner}, which are applied from the right. GMRES restarts after 30 iterations, with a maximum of 120 iterations, and a relative tolerance is $10^{-12}$. 
All AMG subsolvers considered use the strong threshold of 0.7 and a single V-cycle. 
The restriction operators for mechanics AMG are built by BoomerAMG based on a submatrix of a single displacement component, and the same restriction operators are applied to each component \cite{hmg}. Two iterations of point-block smoothers from PETSc are applied on each AMG levels: SOR for 2D and ILU(0) for 3D. The direct solver is used in the coarsest level.
The pressure AMG solver in the CPR method and the coupled pressure-temperature AMG subsolver in the System-AMG method are fully are using BoomerAMG with a default symmetric-SOR smoothers. PETSc implementation of point-block ILU(0) is used as the local preconditioner $M_L$, the second part of the CPR method. Before applying it, the pressure and temperature variables are reordered into an interleaved arrangement.

To ensure the validity of the Schur complement approximations in \Cref{sec:eliminating_intf_flow,sec:eliminating_force_balance}, we additionally run experiments 
where the inner subsolvers are accelerated by GMRES iterations with relative tolerance of $10^{-5}$, again using right preconditioning.
As this leads to a non-stationary preconditioner, the outer system is solved with 
the flexible GMRES (FGMRES) algorithm. %
This solver uses CPR as a pressure-temperature preconditioner. 
Though this algorithm results in a much larger solution time, the subproblem solvers are close to being exact. 
Therefore, the only sources of inexactness in the preconditioner are the Schur complement approximations, and the FGMRES iteration count under grid refinement serves as a measure of the scalability of these approximations.

The code for the runscripts of the numerical experiments can be found on GitHub\footnote{\href{https://github.com/Yuriyzabegaev/FTHM-Solver/}{https://github.com/Yuriyzabegaev/FTHM-Solver/}}, along with the Docker image \cite{zabegaev2025runscripts}. It is based on PorePy, the mixed-dimensional porous media simulation framework \cite{porepy2024}, implemented in Python. The computational simplex grids are generated using Gmsh \cite{gmsh}. The preconditioner algorithms are implementing in PETSc \cite{petsc4py,petsc_user_manual}, interfacing BoomerAMG from HYPRE \cite{hypre}.

\begin{table}[ht]
    \centering
\begin{tabular}{llll}
\toprule
\textbf{Symbol} & \textbf{Value} & \textbf{Unit} & \textbf{Parameter Name} \\
\midrule
$\lambda_{\text{Lame}}$ & $1.2 \cdot 10^{10}$ & Pa & \LameName first parameter \\
$G$ & $1.2 \cdot 10^{10}$ & Pa & \LameName second parameter \\
$\alpha$ & $0.47$ & - & Biot coefficient \\
$K_m$ & $1.0 \cdot 10^{-13}$ & m$^2$ & Matrix isotropic permeability \\
$K_j$ & $1.0 \cdot 10^{-4}$ & m$^2$ & Interface normal permeability \\
$\gamma$ & $4.559 \cdot 10^{-10}$ & Pa$^{-1}$ & Fluid compressibility \\
$\mu^f$ & $1.002 \cdot 10^{-3}$ & Pa s & Fluid viscosity \\
$F$ & $0.577$ & - & Friction coefficient \\
$\theta$ & $5.0$ & $^\circ$ & Shear dilation angle \\
$c^f$ & $4182.0$ & J kg$^{-1}$K$^{-1}$ & Fluid specific heat capacity \\
$c^s$ & $720.7$ & J kg$^{-1}$K$^{-1}$ & Solid specific heat capacity \\
$\kappa^f$ & $0.5975$ & kg m$^{-1}$ & Fluid heat conductivity \\
$\kappa^s$ & $0.1$ & kg m$^{-1}$ & Solid heat conductivity \\
$\beta^f$ & $2.086 \cdot 10^{-4}$ & Pa s & Volumetric fluid thermal expansion \\
$\beta^s$ & $9.66 \cdot 10^{-6}$ & Pa s & Volumetric solid thermal expansion \\
$p_{\text{ref}}$ & $3.5 \cdot 10^7$ & Pa & Reference pressure \\
$\rho_{\text{ref}}^f$ & $998.2$ & kg & Reference fluid density \\
$a_{\text{ref}}$ & $1.0 \cdot 10^{-3}$ & m & Reference aperture \\
$\varphi_{\text{ref}}$ & $1.3 \cdot 10^{-2}$ & - & Reference porosity \\
$T_{\text{ref}}$ & $393$ & K & Reference temperature \\
$\rho_{\text{ref}}^s$ & $2683.0$ & kg & Reference solid density \\
\bottomrule
\end{tabular}

    \caption{The material constants used in the numerical experiments. $K_m$ represents the matrix isotropic permeability. The remaining symbols are formally defined in \Cref{appendix:model_description}.}
    \label{tab:material_params}
\end{table}

\subsection{2D model with fractures}

\begin{figure}
    \centering
    \includegraphics[width=1\linewidth]{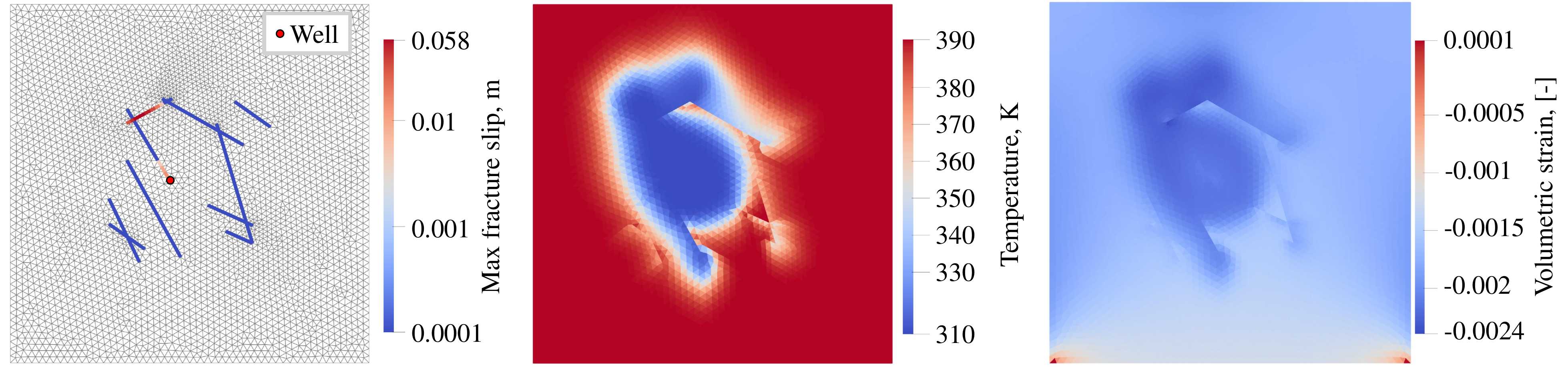}
    \caption{State of the 2D model using medium coarse grid with 3e+4 degrees of freedom at $t \approx 375$ days. Left: The computational grid with bold 1D fractures and maximum fracture slip. The red dot marks the injection cell.
    Center: Temperature distribution. Right: Volumetric strain distribution.}
    \label{fig:model_2d}
\end{figure}

The first experiment serves as a proof-of-concept for the preconditioners applied to a 2D problem. 
The simulation domain is a 2D square with side lengths $2$ km. The domain is compressed by a normal force on the north side at $8 \cdot 10^7$ Pa and $9.6 \cdot 10^7$ Pa on the west and the east sides. The south side is fixed.

The fracture geometry follows Benchmark 3 of \cite{FLEMISCH2018239}, consisting of ten line fractures and six intersections forming three disconnected fracture clusters. The benchmark geometry is scaled to fit within a $1$-km square and positioned at the domain center, ensuring that each fracture remains at least $500$ m away from the domain boundary. The injection well is located in the fracture cell closest to the domain center, constantly injecting the cold fluid with the rate of $10 \text{ kg s\textsuperscript{-1}}$.  The simulated time is 500 days. 
The setup geometry, positions of slipping fractures, and the temperature distribution are illustrated in \Cref{fig:model_2d}. 

To confirm that the fracture deformation is impacted by thermal effects, and not only the injection pressure, we also ran a simulation with the same setup, but with the injection temperature set to the reference temperature. This lead to no sliding of fractures, while in the simulations reported below, approximately 10\% of the fracture cells slid during the simulation.

\subsubsection{Scaling under grid refinement}
\begin{figure}
    \centering
    \includegraphics[width=1\linewidth]{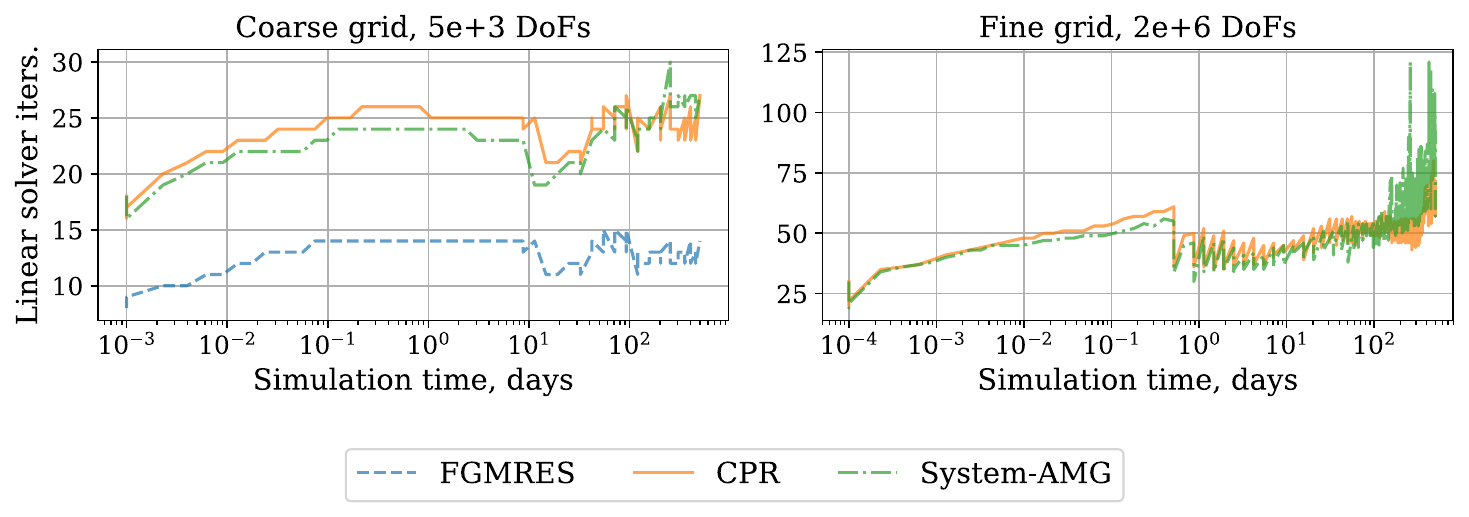}
    \caption{Number of GMRES iterations in the 2D simulation, plotted against the simulation time, for the coarsest (left) and finest (right) grids. Results with CPR and System-AMG are shown for both grids. The FGMRES solver, applied only on the coarsest grid, has subproblem solvers that are close to being exact.
    The observed increase in GMRES iteration counts is due to the adaptive time-stepping strategy, wherein the simulation starts with a time step of 1e-3 seconds and progressively grows to as large as 3 years.
    }
    \label{fig:gmres_iters_2d}
\end{figure}

We first tested the scaling behavior of the preconditioners when applied to a sequence of increasingly finer grids.
\Cref{fig:gmres_iters_2d} shows how the number of GMRES iterations vary through the simulation on the coarsest and finest grid for the CPR and System-AMG preconditioners.
Though there are minor differences between the preconditioners, the two preconditioners show very similar performance through most of the simulations. The only notable exception is that System-AMG seems to struggle towards the end of the simulation on the finest grid. However, when considering the average number of iterations throughout the simulation, presented in \Cref{tab:exp2_gmres_iters}, the difference is minor, and there is no systematic trend in the results. For both preconditioners there is a modest increase in the iteration count, which roughly doubles as the number of degrees of freedom is increased with almost three orders of magnitude.

\Cref{fig:gmres_iters_2d} and \Cref{tab:exp2_gmres_iters} also show the number of FGMRES iterations, representing a preconditioner with close to exact subsolvers. Due to the high computational cost of this preconditioner, we did not apply these on the two finest grids.
The number of FGMRES iterations can be seen to be stable under grid refinement and show minor variations throughout a simulation. This indicates that the Schur complement approximation is scalable, and that the increase in iteration count is caused by the subsolvers for individual components of the Schur complement.

\begin{table}[t]
    \centering
\begin{tabular}{lrrrrrr}
\toprule
Total DoFs & 5e+3 & 9e+3 & 3e+4 & 6e+5 & 1e+6 & 2e+6 \\
\midrule
CPR & 23.9 & 24.8 & 28.1 & 43.5 & 48.4 & 50.6 \\
System-AMG & 23.4 & 24.4 & 30.9 & 46.5 & 45.9 & 55.6 \\
FGMRES & 12.7 & 12.3 & 11.6 & 10.8 & - & - \\
\bottomrule
\end{tabular}
    \vspace{5pt}
    \caption{Average number of linear solver iterations with the variants of the preconditioner in the 2D experiment with fractures. FGMRES denotes the application of the corresponding linear solver with the inner GMRES subsolvers and relative tolerance of 1e-5. Experiments marked ``-'' are skipped.}
    \label{tab:exp2_gmres_iters}
\end{table}

\subsubsection{Dependency on frictional coefficient}
To further probe the preconditioners under varying regimes of fracture deformation, we next considered a series of simulations where the friction coefficient $F$ is varied from a low value of 0.1 to a high value of 0.9. 
Low values of $F$ give low resistance to slip, while with higher values, the fractures are more likely to be in a sticking regime.
In practice, the friction coefficient for rocks is commonly taken to be between 0.5 and 0.7, e.g. \cite{Jaeger2005}.
This experiment is ran with the grid refinement fixed at level four, corresponding to 6e+5 degrees of freedom.

\Cref{tab:friction} presents the maximum number of fracture cells sliding simultaneously and the corresponding linear solver iterations. The results indicate that amplified sliding increases problem complexity for both the CPR and System-AMG preconditioner variants, but the increase in iteration count is limited even for the lowest value of $F$, when almost all fracture cells undergo sliding. 
Though System-AMG performs best for low values of $F$, there are no marked differences for more realistic friction values.

\subsubsection{Dependency on Peclet number}
As demonstrated in \cite{roy_block_2019, cremon_constrained_2024}, AMG-based preconditioners for temperature transport have been shown to outperform CPR in simulations of processes spanning hundreds of years, where the primary driving force is diffusion rather than injection. 
The limited domain size and simulation time and the relatively high injection rate makes the current setup injection dominated, as is reflected in the performance of CPR.

To test the preconditioners under different regimes, we respectively increased and decreased the thermal conductivity of the solid with two orders of magnitude, and reran the simulation at the fourth refinement level, corresponding to 6e+5 degrees of freedom.
As can be seen from the averages of linear iterations, shown in \Cref{tab:peclet}, there is no significant decline in CPR performance.
This suggests that CPR is suitable for the current injection-driven setup across a broad range of material coefficients, making it the preferable preconditioner, as each iteration of CPR requires less computation compared to System-AMG.
However, if modeling long-term processes becomes necessary, and CPR's limitations are expected, System-AMG remains a viable option.
We will return to this question in the 3D experiment.

\begin{table}[ht]
    \centering
    \begin{minipage}{0.49\textwidth}
        \centering
        \begin{tabular}{lrrrrr}
\toprule
Friction $F$ & 0.1 & 0.2 & 0.3 & 0.577 & 0.9 \\
\midrule
Sliding, \% & 96 & 54 & 24 & 9 & 4 \\
CPR & 66.8 & 55.1 & 47.5 & 43.5 & 46.0 \\
System-AMG & 56.6 & 50.2 & 43.3 & 46.5 & 47.9 \\
\bottomrule
        \end{tabular}
        \vspace{5pt}
        \captionof{table}{Values of the friction coefficient $F$, percentage of fracture cells at the sliding state simultaneously, and average number of GMRES iterations with the variants of the preconditioner.}
        \label{tab:friction}
    \end{minipage}
    \hfill
    \begin{minipage}{0.49\textwidth}
        \centering
        \begin{tabular}{lrrr}
        \toprule
        Diffusion $\kappa^s$ & 1e-3 & 1e-1 & 1e+1 \\
        \midrule
        Pe & 1e+6 & 7e+4 & 1e+3 \\
        CPR & 43.6 & 43.5 & 42.4 \\
        System-AMG & 45.6 & 46.5 & 47.0 \\
        \bottomrule
        \end{tabular}
        \vspace{5pt}
    \caption{Various values of the diffusion coefficient $\kappa^s$ and average number of GMRES iterations with the variants of the preconditioner. `Pe' denotes the \PecletName{} number, computed as the ratio of the maximum convection and diffusion fluxes in the domain and averaged over time.}
    \label{tab:peclet}
    \end{minipage}

\end{table}

\begin{figure}[t]
    \centering
    \includegraphics[width=1\linewidth]{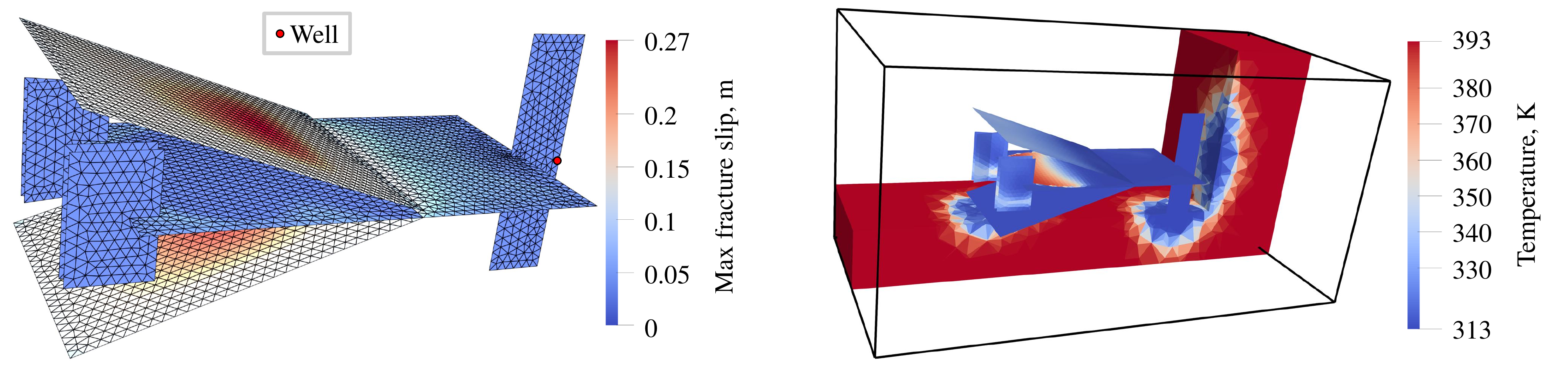}
    \caption{State of the 3D model using grid with 6e+5 degrees of freedom at $t \approx 30$ years. Left: Computational grid of the fracture subdomains and maximum fracture slip. The upper-left fracture is disconnected from the remainder of the fracture network with a small gap. Right: Temperature distribution within fractures and in the slice of the porous medium.}
    \label{fig:model_3d}
\end{figure}

\begin{figure}
    \centering
    \includegraphics[width=1\linewidth]{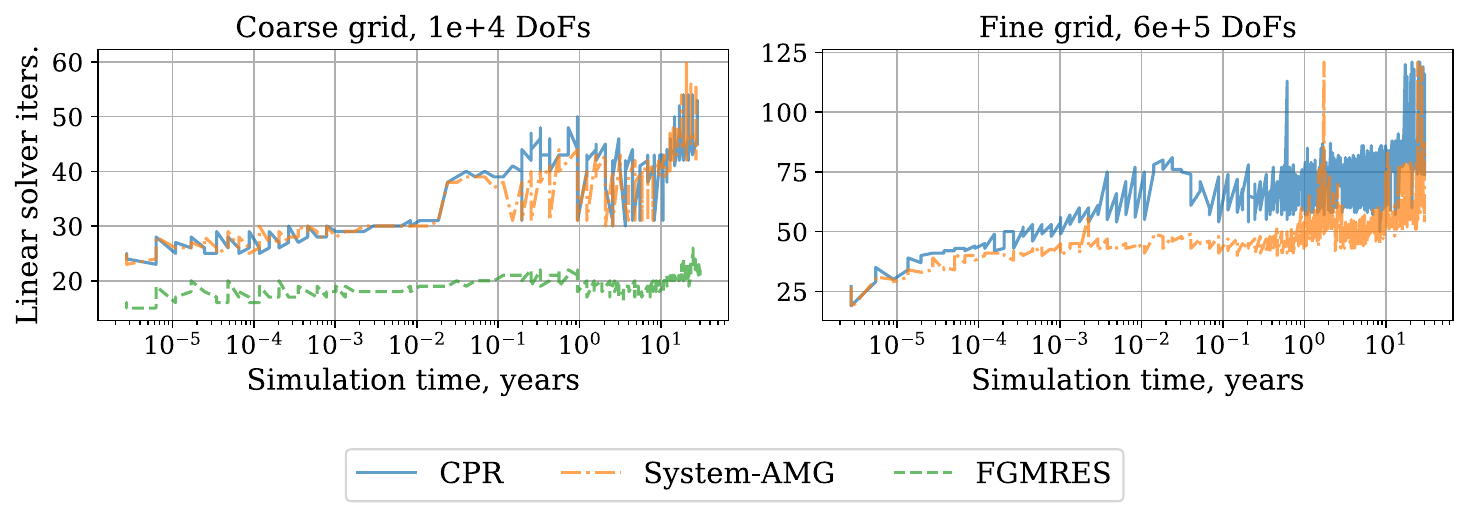}
    \caption{Number of GMRES iterations in the 3D simulation, plotted against the simulation time, for the coarsest (left) and finest (right) grids. Results with CPR and System-AMG are shown for both grids. The FGMRES solver, applied only on the coarsest grid, has subproblem solvers that are close to being exact. The observed increase in GMRES iteration counts is due to the adaptive time-stepping strategy, as the simulation starts with a time step of 1e-3 seconds and progressively grows to as large as 3 years.}
    \label{fig:gmres_iters_3d}
\end{figure}

\subsection{3D model with fractures}

The second model consists of a 3D box domain measuring $1$ km along the x and z axes and $2$ km along the y axis. The fracture geometry follows Case 3 of \cite{BERRE2021103759}, consisting of eight plane fractures and seven line intersections. All fractures, except for one, are connected, while the disconnected fracture is positioned close to the others. The fracture network is scaled to fit within a box that is half the size of the model and is positioned at the center of the domain, so the distance between the fracture and the domain boundary is at least $500$ m. Again we consider a sequence of increasingly finer grids.

The deformation is fixed at $y = y_\text{max}$, while the other boundaries are subjected to a constant normal compressive force: $8\cdot 10^7$ Pa along $y = 0$; $1.2 \cdot 10^8$ Pa along $x = 0$ and $x = x_\text{max}$; $4 \cdot 10^7$ Pa along $z = 0$ and $z = z_\text{max}$. The inlet is connected to the fracture cell closest to the $y = 0$ boundary. The simulated time is 30 years, and the injection rate is 100 kg s\textsuperscript{-1}. 

The fractures positions and the sharp angles at intersections make this geometry challenging for meshing and inherently difficult for the linear solver, so we use the ILU(0) as a smoother for the mechanical AMG instead of SOR. 

\begin{figure}[ht]
    \centering
    \begin{minipage}{0.45\textwidth}
        \centering
        \includegraphics[width=\linewidth]{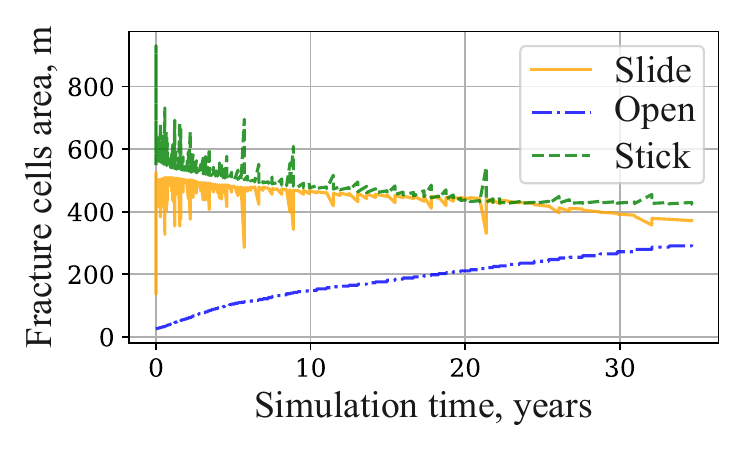} %
        \caption{Evolution of fracture cell states in the 3D simulation. Injection process initiated sliding, and fractures increasingly open due to the fluid pressure.}
        \label{fig:fracture_area}
    \end{minipage}
    \hfill
    \begin{minipage}{0.45\textwidth}
\begin{tabular}{lrrrr}
\toprule
Total DoFs & 1e+4 & 3e+4 & 9e+4 & 6e+5 \\
\midrule
CPR & 38.0 & 38.0 & 42.6 & 70.9 \\
SAMG & 36.8 & 33.6 & 37.2 & 51.0 \\
FGMRES & 19.5 & 20.7 & 18.7 & - \\
\bottomrule
\end{tabular}
    \vspace{5pt}
    \captionof{table}{Average number of linear solver iterations with the variants of the preconditioner in the 3D experiment with fractures. FGMRES denotes the application of the corresponding linear solver with the inner GMRES subsolvers and relative tolerance of 1e-3. Experiment marked ``-'' is skipped.}
    \label{tab:exp3_gmres_iters}
    \end{minipage}
\end{figure}

\Cref{fig:model_3d} shows the setup geometry together with the temperature distribution at the end of the simulation.
Compared to the 2D simulation, the injection rate is lower compared to the domain volume. This is reflected in the cold water front caused by the fluid injection having intruded a relatively short distance from the fracture network into the host rock, even at the end of the simulation.
This lower injection rate implies that diffusion may play a comparatively larger part as a mechanism for energy transport, compared to the 2D case. The distribution of fracture states over simulation time is shown in \Cref{fig:fracture_area}, highlighting that the simulation includes states where up to 55\% of fracture cells were sliding and up to 30\% were opening simultaneously.

Consider first the number of GMRES iterations as a function of simulation time, shown in \Cref{fig:gmres_iters_3d}.
For the coarsest grid, CPR and System-AMG perform similar. For both methods, there is a notable increase in iteration count towards the end of the simulation, however, this coincides with increasing time step sizes due to the adaptive time-stepping scheme.
On the finest grid, there is however a notable difference between the methods, with System-AMG requiring substantially fewer iterations than CPR. We believe this is related to the larger role played by diffusion in this simulation, in particular after the pulse of elevated pressure has propagated through the fracture network in the initial part of the simulation.
When considering the average numbers of iterations for all grid refinement, see \Cref{tab:exp3_gmres_iters}, it is clear that System-AMG shows a better scaling behavior than CPR does with a $1.4$ times increase in GMRES iteration number against $1.9$ times.
Nevertheless, as System-AMG also has a higher cost of setup and application, CPR remains a viable option also for the current setup.
As shown in \Cref{fig:gmres_iters_3d} and \Cref{tab:exp3_gmres_iters}, the number of FGMRES iterations remains stable, indicating that the good scalability of the Schur complement approximations is preserved in 3D.

\section{Conclusion}
\label{sec:results}
This study presents a linear preconditioner for the thermo-poromechanics problem in fractured porous media, accounting for fractures that can undergo frictional sliding. The preconditioner is based on the fixed-stress based approximation, extended for the porous media and fracture subdomains, to decouple the momentum balance, while the coupled mass and energy balance is treated by either one of the two variations: CPR and System-AMG. The contact mechanics subproblem is decoupled with help of the linear transformation, which ensures that the corresponding submatrix is not singular. 

Both variations of the preconditioner are extensively tested by the numerical experiments. The preconditioner variant with System-AMG demonstrates better scaling for the complex 3D model under computational grid refinement, but
the preconditioner variant with CPR still remains a viable choice due to its lower computational cost.
Additionally, the simulations in the grid refinement experiment have passed through various fracture states including simultaneous sliding of the majority of the fracture cells,
indicating that the preconditioner works well regardless of the fracture state.
Finally, the algorithms showed robust performance across various combinations of realistic model parameters: Varying the \PecletName number 
lead to negligible changes in GMRES iteration number for both variations of the algorithm, while changing the fracture friction coefficient 
caused only a slight change for both System-AMG and CPR variants of the preconditioner.

\section*{Acknowledgment}
This project has received funding from the VISTA program, The Norwegian Academy of Science and Letters and Equinor, from the Research Council of Norway (grant 308733), and from the European Research Council (ERC) under the European Union’s Horizon 2020 research and innovation programme (grant agreement No 101002507).

\begin{appendices}
\crefalias{section}{appendix}

\section{Model description}
\label{appendix:model_description}

Following the presentation of \cite{porepy2024}, this appendix presents the mathematical model for mixed-dimensional flow, deformation, and energy transfer in a poroelastic medium, incorporating frictional fracture deformation and governed by contact poromechanics, previously introduced in \cite{stefansson_fully_2021,berre2019flow,porepy2024,berge_finite_2020,porepy2021}. 
The appendix consists of the following parts: \Cref{sec:mix_dim_geo} revisits the mixed-dimensional model geometry, while \Cref{sec:conservation_laws} presents the laws for the momentum conservation and the fluid mass and energy conservation. Next, the contact mechanics relations, including the Coulomb friction law and the non-penetration condition, are described in \Cref{sec:model_description_contact_mech}, and the constitutive laws are listed in \Cref{sec:constitutive_laws}. Finally, the initial and the boundary conditions are defined in \Cref{sec:boundary_conditions}

\subsection{Geometry of a mixed-dimensional model}
\label{sec:mix_dim_geo}
The model domain is represented as a collection of subdomains $\Omega_i$ with different varying $d_i$, where $d_i \in \left\{0, ..., D \right\}$ and $D \in \left\{ 2, 3 \right\}$. 
When $D = 3$, the model can incorporate several subdomains of varying dimensions: the porous medium, represented by a 3D subdomain; the fractures, modeled as 2D subdomains; the intersections of fractures, which are treated as 1D subdomains; and the intersections of intersections, represented as 0D subdomains.
The width of a fracture subdomain $\Omega_i$ is characterized by its aperture $a_i$. To represent the spatial extension of $\Omega_i$ with dimension $d_i < D$ orthogonal to $\Omega_i$, we define the specific volume $\nu_i = a_i ^ {D - d_i}$. For fracture subdomains, $\nu_i$ equals to the aperture $a_i$, and for fracture intersections, it is equal to $a_i^2$.
For the porous medium subdomain with $d_i=D$, we define $\nu_i=1$ for consistency.

The coupling between the neighboring subdomains of co-dimension 1 is achieved through an interface $\Gamma_j$. We denote the higher-dimensional and lower-dimensional neighboring subdomains of an interface with subscripts $h$ and $l$ when applicable. The boundary of $\Omega_i$ is denoted as $\partial \Omega_i$, while the internal portion of $\partial \Omega_i$ that geometrically coincides with the interface $\Gamma_j$ is written as $\partial_j \Omega_i \subseteq \partial \Omega_i$. The projection operator that maps relevant quantities from subdomain $\Omega_i$ to interface $\Gamma_j$ is denoted by $\Pi^i_j$ and the reverse operator is denoted by $\Xi^i_j$.

Subscripts are used to specify the subdomain or interface where the quantities are defined, while superscripts $f$ and $s$ indicate fluid and solid quantities, respectively. In cases where the context is clear, subscripts and superscripts are omitted for simplicity.

\subsection{Conservation laws}
\label{sec:conservation_laws}

The following subsection presents the fluid mass, energy and momentum conservation laws for the relevant subdomains.
The fluid mass conservation equation for the subdomains $\Omega_i$ of dimension $d_i \in \left\{ 0, ..., D \right\}$:

\begin{equation}
    \label{eq:mass_conservation}
    \frac{\partial}{\partial t} \left( \nu_i \rho_i^f \varphi_i \right) - \nabla \cdot \left( \nu_i \rho_i^f \mathbf{v}_i \right) -\sum_{j \in \hat{S}_i} \Xi^i_j \left( \nu_j \rho_j^f v_j \right) = \psi_i.
\end{equation}
\noindent
Here, $\rho_i^f$ and $\rho_j^f$ represent the fluid densities defined within a subdomain and on an interface, respectively; $\varphi_i$ denotes the porosity; $\mathbf{v}_i$ and $v_j$ are the volumetric fluid fluxes within a subdomain and on an interface, respectively; $\psi_i$ represents a source or sink of fluid mass; $\nu_j \coloneq \Pi_i^h \nu_h$ is the specific volume of the interface; and the set $\hat{S}_i$ includes all interfaces connecting to higher-dimensional neighbors of $\Omega_i$. Fluxes to lower-dimensional neighbors are treated as internal boundary conditions.
The second term vanishes when $d_i = 0$, as 0D domains do not permit internal flux, 
while the third term is absent when $d_i = D$ as these subdomains lack higher-dimensional neighbors.

The energy conservation equation for the subdomains $\Omega_i$ of dimension $d_i \in \left\{ 0, ..., D \right\}$:
\begin{equation}
    \label{eq:energy_conservation}
    \frac{\partial}{\partial t} \left(
    \nu_i \varphi_i U^f_i + \nu_i (1 - \varphi_i) U^s_i 
    \right)
    +
    \nabla \cdot 
    (\bm{w}_i + \bm{q}_i)
    + \sum_{j \in \hat{S}_i} \Xi^i_j (w_j + q_j)
    = \psi^T_i.
\end{equation}
Here,
$U^f_i$ and $U^s_i$ are the fluid and solid internal energies, $\bm{w}_i$ and $\bm{q}_i$ are the convective and diffusive energy flux terms within a subdomain, and $w_j$ and $q_j$ are the inter-dimensional convective and diffusive flux terms, respectively. $\psi^T_i$ denotes the energy source term.

The quasi-static momentum conservation equation is applied to the porous medium subdomain ($d_i = D$):
\begin{equation}
    \label{eq:momentum_conservation}
    -\nabla \cdot \sigma_i  = \mathbf{F}_i.
\end{equation}
Here, $\sigma_i$ denotes the total stress tensor, while $\mathbf{F}_i$ denotes the body forces.

\subsection{Contact mechanics relations}
\label{sec:model_description_contact_mech}

This subsection considers a pair of the porous medium $\Omega_h$ and a fracture $\Omega_l$ of dimensions $d_h = D$ and $d_l = D - 1$, respectively. Two interfaces are denoted on each side of the fracture: $\Gamma_j$ and $\Gamma_k$. A normal vector on a fracture $\mathbf{n}_l$ is defined to coincide with $\mathbf{n}_h$ on the $j$-side of the fracture. We introduce the contact traction variable $\lambdafull_l$, defined in $\Omega_l$ according to the direction of $\mathbf{n}_l$. For a general 3D vector $\mathbf{\iota}$, the tangential and normal components are defined as follows:
\begin{equation}
    \mathbf{\iota}_\mathrm{n} = \mathbf{\iota} \cdot \mathbf{n}_l, \quad \mathbf{\iota}_\tau = \mathbf{\iota} - \mathbf{\iota}_\mathrm{n}.
\end{equation}
We also introduce the jump in interface displacements across $\Omega_l$:
\begin{equation}
    \ujump = \Xi^l_k \disp_k - \Xi^l_j \disp_j.
\end{equation}
According to Newton's third law, there is a balance between the contact traction and the total traction on each fracture surface:
\begin{equation}
\begin{aligned}
\label{eq:intf_force_balance}
        \Pi^h_j \sigma_h \cdot \mathbf{n}_h & = \Pi^l_j \left( \lambdafull_l - p_l \mathbf{I} \cdot \mathbf{n}_l \right),\\
            -\Pi^h_k \sigma_h \cdot \mathbf{n}_h &= \Pi^l_k \left( \lambdafull_l - p_l \mathbf{I} \cdot \mathbf{n}_l \right).
\end{aligned}
\end{equation}
\noindent
Here, $\mathbf{I}$ denotes the identity matrix.
We omit the subscript $l$ for readability and enforce the non-penetration condition for the fracture surfaces:
\begin{equation}
\label{ineq:contact_normal_ineq}
   \ujumpn - g \geq 0;
   \quad
    \lambdan \leq 0;
    \quad
    \lambdan \left( \ujumpn - g \right) = 0.
\end{equation}
Here, $g$ defines the normal gap between the fracture's surfaces when they are in mechanical contact. 
The second inequality indicates that compressive normal contact traction is associated with a negative $\lambdan$.
The friction bound, $b$, is defined according to the Coulomb friction law with a constant friction coefficient $F$:
\begin{equation}
\label{eq:friction_bound}
    b = -F \lambdan.
\end{equation}
The velocity of the tangential displacement jump in case of sliding is denoted $\ujumpt$, while the friction model is given by:
\begin{equation}
\begin{aligned}
\label{ineq:contact_tangential_ineq}
    \norm{\lambdat} &\leq b, \\
    \norm{\lambdat} &< b \quad \rightarrow \quad \ujumpt = 0, \\
    \norm{\lambdat} &= b \quad \rightarrow \quad \exists \zeta \in \mathbb{R}^+ : \ujumpt = \zeta \lambdat.
\end{aligned}
\end{equation}
The three relations express that: (i) tangential stresses are bounded, (ii) tangential deformation only occurs when the bound is reached, and (iii) tangential stresses and deformation increments are co-directed.
The contact mechanics inequalities determine that every point on the fracture is categorized into one of three distinct states:
\begin{itemize}
    \item Sticking -- the fracture surfaces are in mechanical contact and remain stationary relative to one another; the traction $\lambdafull$ is nonzero.
    \item Sliding -- the fracture sides are in mechanical contact with the tangential contact force $\lambdat$ having reached the friction bound, and the fracture surfaces are moving relative to each other.
    \item Open -- the fracture surfaces are separated, with no mechanical contact is occurring between them; the traction $\lambdafull$ is zero.
\end{itemize}
The inequalities, \cref{ineq:contact_normal_ineq,ineq:contact_tangential_ineq}, can be reformulated as equalities using complementarity functions, in accordance with the augmented Lagrangian formulation, as detailed in \cite{HUEBER20053147, stefansson_fully_2021, berge_finite_2020}:
\begin{subequations}
\label{eq:contact_equalities}
    \begin{equation}
        \lambdan + \text{ max}\left\{ 0, -\lambdan - c (\ujumpn - g)  \right\} = 0;
    \end{equation}
    \begin{equation}
        \left( -\lambdat \text{max}\left\{ b, \norm{\lambdat + c\ujumpt} \right\} + \text{max}\left\{ b, 0 \right\} \left(\lambdat + c \ujumpt \right) \right) \left(1 - \chi \right) + \chi \lambdat = 0.
    \end{equation}
\end{subequations}
Here, $c>0$ is a numerical parameter, and $\chi$ is a characteristic function defined as: $\chi = 1 \text{ if } |b| < \varepsilon$ and $\chi = 0 \text{ otherwise}$. Here,
$\varepsilon$ is a numerical parameter that treats the fracture as open even for small positive values of the friction bound. This approach helps reduce the sensitivity of the nonlinear solver to numerical errors, thereby improving convergence.

\subsection{Constitutive laws}
\label{sec:constitutive_laws}

The constitutive relations outlined below are primarily derived from \cite{coussy2004poromechanics}. The volumetric fluid flux is defined according to Darcy's law:
\begin{equation}
\label{eq:darcy}
    \mathbf{v}_i = - \frac{\perm_i}{\visc} \left( \nabla p_i - \rho^f_i \gravity \right).
\end{equation}
Here, $\visc$ denotes fluid viscosity, and $\gravity$ is gravity. The permeability tensor $\perm_i$ is assumed to be constant in the matrix, while the fracture subdomains permeability is given by \cite{zimmerman1996hydraulic}:
\begin{equation}
\label{eq:fracture_permeability}
    \perm_i = \frac{a_i^2}{12}\mathbf{I}, \quad d_i = D - 1,
\end{equation}
where $a_i$ depends on $\ujump$ as detailed in eq. \eqref{eq:aperture}. The permeability of intersections is computed as the average of the permeabilities in the intersecting fractures.
While the solid density is constant, the fluid density is given by (we omit the subscript $i$):
\begin{equation}
\label{eq:fluid_density}
    \rho^f = \rhoRef \exp{\left[ \compres (p - p_\text{ref}) + \beta^f (T - T_\text{ref})\right] },
\end{equation}
where $\compres$ denotes compressibility, $\beta^f$ denotes fluid thermal expansion, $\rhoRef$, $p_\text{ref}$ and $T_\text{ref}$ denote the reference fluid density, pressure and temperature, respectively. The viscosity is constant.

The internal energy of the solid and the fluid are:
\begin{equation}
    U^s = \rho^s c^s T; 
    \quad 
    U^f = \rho^f (\specHeatCapacity T - p)
    ,
\end{equation}
where $c^s$ and $\specHeatCapacity$ are the specific heat capacities of the solid and the fluid, respectively. The convective and the diffusive energy fluxes is given by:
\begin{equation}
    \bm{w}_i = \nu_i \rho_i^f \specHeatCapacity  T_i \mathbf{v}_i;
    \quad
    \bm{q}_i = \kappa^e_i \nabla T_i,
\end{equation}
where the effective heat conductivity is given by $\kappa^e_i = \varphi_i \kappa^f + (1 - \varphi_i) \kappa^s$ for the fluid and solid heat conductivities $\kappa^s$ and $\kappa^f$, respectively. 
The thermal source term equals the internal energy of the fluid of the volume source and sink terms, i.e. $\psi_i^T = \rho_i^f \specHeatCapacity T_i \psi_i$.

The volumetric interface flux is proportional to the pressure difference across $\Gamma_j$ via a Darcy-type law \cite{Martin2005}:
\begin{equation}
\label{eq:interface_flow}
    v_j = - \frac{\mathcal{K}_j}{\visc_{j}} \left[ \frac{2}{\Pi^l_j a_l} \left( \Pi^l_j p_l - \Pi^h_j p_h \right) -\mathbf{g} \cdot \mathbf{n}_h \rho_j^f \right].
\end{equation}
Here, the interface permeability and viscosity $\mathcal{K}_j$ and $\visc_{j}$ and inherited from the lower-dimensional neighbor subdomain, and the multiplier $\frac{2}{\Pi^l_j a_l}$ represents half the normal distance across the fracture. An advective quantity defined on interfaces $\xi_j$, which represents either $\rho_j$ or $\visc_{j}$, is treated by an inter-dimensional upwinding based on $v_j$:
\begin{equation}
    \xi_j = \begin{cases}
        \Pi^h_j \xi_h \quad &\text{if } v_j > 0 \\
        \Pi^l_j \xi_l \quad &\text{otherwise}.
    \end{cases}
\end{equation}
The advective interface energy flux is defined according to the upstream direction of the interface fluid flux:
\begin{equation}
\label{eq:interface_energy_advective}
    w_j = \begin{cases}
        v_j \Pi^h_j \rho^f_h \specHeatCapacity_{,h} T_h \quad &\text{if } v_j > 0 \\
        v_j \Pi^l_j \rho^f_l \specHeatCapacity_{,l} T_l \quad &\text{otherwise}.
    \end{cases}
\end{equation}
The Fourier-type conductive interface flux is:
\begin{equation}
\label{eq:interface_energy_diffusive}
    q_j = -\kappa_j \frac{2}{\Pi^l_j a_l} (\Pi^l_j T_l - \Pi^h_j T_h),
\end{equation}
with the normal heat conductivity modeled as $\kappa_j = \Pi^l_j \kappa^f_l$, since it originates from the dimension reduction of a fluid-filled domain.

The total thermo-poromechanical stress tensor is described by Hooke's law and incorporates fluid pressure and temperature effects (with the subscript $i$ omitted):
\begin{equation}
\label{eq:poromech_stress}
    \sigma = \sigma_\text{ref} + G \left( \nabla \disp + \nabla \disp^T \right) 
    + 
    \Lame \text{tr}(\nabla \disp) \mathbf{I} 
    - 
    \alpha \left( p - p_\text{ref} \right) \mathbf{I}
    -
    \beta^s (2 G + 3 \Lame) (T - T_\text{ref}) \mathbf{I}
    ,
\end{equation}
where $\sigma_\text{ref}$ is the reference stress configuration, $\alpha$ is the Biot coefficient, $G$ and $\Lame$ are the Lame coefficients, $\text{tr}(\cdot)$ denotes the trace of a matrix, $\beta^s$ is the volumetric solid thermal expansion.
The matrix porosity depends on pressure, temperature and displacement according to \cite{coussy2004poromechanics}:
\begin{equation}
\label{eq:porosity}
    \varphi = \varphi_\text{ref} + \alpha \nabla \cdot \mathbf{u} + \frac{(\alpha - \varphi_\text{ref})(1 - \alpha)}{\Lame + \frac{2}{3} G} \left( p - p_\text{ref} \right)
    - {(\alpha - \varphi_\text{ref})} \beta^s (T - T_\text{ref})
    ,
\end{equation}
where $\varphi_\text{ref}$ is the reference porosity.
The porosity in fractures and intersections is set to one.
The gap function is used to model fracture roughness effects:
\begin{equation}
\label{eq:gap_function}
    g = g_\text{ref} + \norm{\disp_\tau} \tan \theta,
\end{equation}
where $g_\text{ref}$ is the steady-state gap and $\theta$ is the shear dilation angle. 

For fracture subdomains, we define the aperture:
\begin{equation}
\label{eq:aperture}
    a = a_\text{ref} + \ujumpn,
\end{equation}
where $a_\text{ref}$ is a residual hydraulic aperture. The aperture of the intersection subdomains is computed as the mean aperture of higher-dimensional neighbors:
\begin{equation}
    a_i = \frac{1}{|\hat{S}_i|}\sum_{j \in \hat{S}_i} \Xi^i_j \Pi^h_j a_h.
\end{equation}

\subsection{Initial and boundary conditions}
\label{sec:boundary_conditions}
To fully specify the system, we set the initial conditions and then impose boundary conditions along both the internal and external perimeters of each subdomain. At the internal interfaces $\partial_j \Omega_i$, we ensure a continuity of the normal mass and energy fluxes, both for advective and diffusive contributions:
\begin{equation}
    \nu_h \rho_h^f \mathbf{v}_h \cdot \vec{n}_h = \Xi^h_j \nu_j \rho^f_j v_j;
    \quad
    \bm{w}_h \cdot \vec{n}_h = \Xi^h_j w_j;
    \quad
    \bm{q}_h \cdot \vec{n}_h = \Xi^h_j q_j.
\end{equation}
At the porous medium subdomain $(d_i = D)$, we enforce the displacement continuity:
\begin{equation}
    \mathbf{u}_i = \Xi^i_j \mathbf{u}_j.
\end{equation}
At the immersed fracture tips, mass and energy fluxes are constrained to zero. On the external boundaries, we impose either Dirichlet or Neumann conditions for the equation of mass, momentum, and energy conservation.

\end{appendices}

\bibliographystyle{unsrt}
\bibliography{references}

\end{document}